\documentstyle[amssymb,amsfonts]{amsart}

\newenvironment{proof}{\noindent {\bf Proof} }{\endprf\par}
\def \endprf{\hfill  {\vrule height6pt width6pt depth0pt}\medskip}
\def\emph#1{{\it #1}}
\def\textbf#1{{\bf #1}}

\newcommand{\bea}{\begin{eqnarray}}
\newcommand{\eea}{\end{eqnarray}}
\def\beaa{\begin{eqnarray*}}
\def\eeaa{\end{eqnarray*}}
\def\c{\cdot}
\def\ba{\begin{array}}
\def\ea{\end{array}}
\def\be#1{\begin{equation} \label{#1}}

\newcommand{\nn}{\nonumber}
\def\medn{\medskip\noindent}
\def\rrrr{{\Bbb R}}
\def\rr{{\bf R}}
\def\rra{\rr^{(1)} }
\def\rrb{\rr^{(2)}}

\def\nn{\nonumber}
\def\stu{{S_{t,u}}}
\def\stau{{S_{\tau,u}}}

\def\err{\mbox{Err}}
\def\DC{D_*}

\def\gg{{\bf g}}

\newcommand{\half}{\frac 12}

\def\f12{\frac 12}
\def\half{\frac 12}

\def\a{\alpha}
\def\b{\beta}
\def\ga{\gamma}
\def\Ga{\Gamma}
\def\de{\delta}

\def\ep{\epsilon}
\def\eps{\epsilon}
\def\la{\lambda}

\def\si{\sigma}
\def\Si{\Sigma}

\def\Om{\Omega}

\def\nab{\nabla}

\def\bP{\overline {P}}

\def\f14{\frac{1}{4}}
\def\dd{{\bf D}}

\newcommand{\les}{\lesssim}

\newcommand{\EE}{{\mathcal E}}

\def\il{\int}

\def\Lb{\underline{L}}

\def\ric{\mbox{ \bf Ric}}

\def\pr{\partial}
\def\f12{\frac 1 2}

\def\trch{\mbox{tr}\chi}

\newcommand{\nabb}{\mbox{$\nabla \mkern-13mu /$\,}}

\def\dtu{D_{t,u}}
\def\dtau{D_{\tau,u}}

\begin{document}
\theoremstyle{plain}
  \newtheorem{theorem}[subsection]{Theorem}
  \newtheorem{conjecture}[subsection]{Conjecture}
  \newtheorem{proposition}[subsection]{Proposition}
  \newtheorem{lemma}[subsection]{Lemma}
  \newtheorem{corollary}[subsection]{Corollary}

\theoremstyle{remark}
  \newtheorem{remark}[subsection]{Remark}
  \newtheorem{remarks}[subsection]{Remarks}

\theoremstyle{definition}
  \newtheorem{definition}[subsection]{Definition}

\include{psfig}

\title[rough Einstein metrics ]{Ricci defects of microlocalized
Einstein metrics}
\author{Sergiu Klainerman}
\address{Department of Mathematics, Princeton University,
 Princeton NJ 08544}
\email{ seri@@math.princeton.edu}

\author{Igor Rodnianski}
\address{Department of Mathematics, Princeton University, 
Princeton NJ 08544}
\email{ irod@@math.princeton.edu}
\subjclass{35J10}
\vspace{-0.3in}
\begin{abstract}
This is the third and last in our series of papers concerning
rough solutions of the Einstein vacuum equations expressed
relative to wave coordinates. In this paper we  prove an important result,
concerning Ricci defects of microlocalized  solutions, stated 
 and used   in the proof the crucial Asymptotics
Theorem in  \cite{Einst2}.

\end{abstract}
\maketitle

\section{Introduction}
This is the third and last in our series of papers concerning
rough solutions of the Einstein vacuum equations expressed
relative to wave coordinates. More precisely we are concerned
with solutions of the Einstein
vacuum equations,
\be{ric1.1}
\rr_{\a\b}(\gg)=0
\end{equation}
expressed\footnote{In wave coordinates
the Einstein equations take the reduced  form 
$$
\gg^{\a\b}\pr_\a\pr_\b  \gg_{\mu\nu}=N_{\mu\nu}(\gg,\pr \gg)
$$
with $N$ quadratic in the first derivatives
 $\pr\gg$ of the metric.} relative to  wave coordinates $x^\a$,
\be{ric1.2}
\square_{\gg} x^\a =\frac{1}{|\gg|}\pr_\mu(\gg^{\mu\nu}|\gg|\pr_\nu)x^\a= 0.
\end{equation}
The  solutions we consider here have a limited degree of differentiability,
we only assume that  in a time slab $[0,T]\times \rrrr^3$ we control the
 the first derivatives of $\gg$ in the  energy norm
$L_t^\infty( H_x^{1+\ga})$, $\ga>0$, as well as in  the mixed
  Strichartz norm $L_t^2(L_x^\infty)$.  More precisely,

\medn
{\bf Metric Hypothesis:}
\be{bootstrap} 
\|\pr\gg\|_{L^\infty_{[0,T]} H^{1+\ga}} +\|\pr \gg\|_{L^2_{[0,T]}
L^\infty_x}\le B_0,
\end{equation}
for some fixed $\ga>0$ arbitrarily small. 

This  condition was introduced in  section 2 of  \cite{Einst1}
as the main  bootstrap assumption in the proof of our main theorem concerning
$H^{2+\ga}$ solutions , $\ga>0$ arbitrarily small, of \eqref{ric1.1}--
\eqref{ric1.2}. 

Microlocalization is  an essential technique  in
 dealing with rough solutions of nonlinear wave equations, see \cite{Einst1}
 and the references therein.
By a  microlocalized  rough Einstein metric, at cut-off parameter
$\la\ge 1$,  we understand, essentially,  the
low frequency part( frequency $<\la$)  of  a given  Einstein 
metric  \eqref{ric1.1}--\eqref{ric1.2}.
To explain this in more details   
we recall below the definition of the 
 Littlewood -Paley   projections,  
\beaa
P_{<\la}&=&\sum_{\mu <\f12 \la}P_{\mu}\\
P_{\mu}f(x)&=&\int e^{i x\c \xi} \chi(\mu^{-1}\xi)\hat{f}(\xi) d\xi
\eeaa 
where $\chi\in {\cal C}_0^\infty(\rrrr)$ supported in $\f12 \le |\xi|\le 2$
and $\sum_{\mu\in 2^{\Bbb Z}}\chi(\mu^{-1}\xi)=1$. The operators $P_\mu$
are the standard Littlewood -Paley {\sl dyadic  projections} corresponding to
 the frequencies $\mu\in 2^{\Bbb Z}$.

 Consider a fixed solution $\gg$  of
\eqref{ric1.1} satisfying the metric hypothesis \eqref{bootstrap}
relative to the fixed  system of  wave coordinates \eqref{ric1.2}.
Consider also a fixed dyadic parameter  $\la\in 2^{\Bbb Z_+}$ and
define the microlocalized  rescaled  metric,
\be{ric1.3}
 H(t,x)=H_{(\la)}(t,x)=(P_{<\la}\gg)(\la^{-1}t, \la^{-1}x)
\end{equation}
Observe that $H_{(\la)}$ is the low frequency part of the 
rescaled metric, i.e. $H_{(\la)}=P_{<1}(G_{(\la)})$ where,
\be{ric1.4}
G_{(\la)}(t,x)=\gg(\la^{-1}t, \la^{-1}x)
\end{equation}

In the rescaled variables we
 restrict ourselves to the slab $[0, t_*]\times \rrrr^3$
with $\,t_*\approx \la^{1-8\ep_0}$ for some small $\ep_0$,
in fact $5\ep_0<\ga$. In this region we define the optical function
$u$ to be the solution of the eikonal equation, 
\be{ric1.5}
H^{\a\b}\pr_\a u\pr_\b u=0
\end{equation}
verifying the initial condition
\be{ric1.6}
u(\Ga_t)=t
\end{equation}
where $\Ga_t$ is the timelike geodesic
passing through  the origin of and orthogonal( with respect to $H$) to
 the initial hypersurface 
$\Si_0$. We denote by $\Si_t$ the spacelike level  hypersurface
generated by the time function $t=x^0$. We denote by $C_u$
the level hypersurfaces of $u$ and by $S_{t,u}$ their intersection
with $\Si_t$. In \cite{Einst2} we show that 
the  the null hypersurfaces $C_u$ form a proper foliation
of the domain $\Om_*={\cal I}^+_{-1}\cap ([0,t_*]\times \rrrr^3)$.
Here ${\cal I}^+_{-1}$ denotes the future domain( domain 
of influence) of the point
$\Ga_{-1}\in \Si_{-1}$.

To each point  $p\in\Om_*$  we associate the canonical null
pair,
\be{ric1.7}
L=T+N,\qquad\qquad \Lb=T-N
\end{equation}
where $T$ is the future unit normal to $\Si_t$
and $N$ is the outward  unit normal to the surface $S_{t,u}$
passing through $p$. Observe that $L$ is proportional to 
the null geodesic generator $L'=-H^{\a\b} \pr_\b u\pr_\a$
of $C_u$.

  A null frame $e_1, e_2, e_3=\Lb, e_4=L$ consists of the
 null pair $L, \Lb$ together with an arbitrary
choice of vectors $(e_A)_{A=1,2}$ tangent to  $S_{t,u}$ such 
that $H(e_A,e_B) =\de_{AB}$. Relative to such a null frame
the metric $H$ has the form,
\be{ric1.8}
H_{34}=-2, \quad H_{33}=H_{44}=H_{3A}=H_{4A}=0,\quad H_{AB}=\de_{AB}.
\end{equation}
The  null components of the  inverse metric are therefore,
\be{ric1.9}
H^{34}=-\f12, \quad H^{33}=H^{44}=H^{3A}=H^{4A}=0,\quad H^{AB}=\de^{AB}.
\end{equation}

 While the  rescaled spacetime metric $G=G_{(\la)}$  verifies the
 Einstein equations
$\rr_{\mu\nu}(G)=0$ this is certainly not true 
for the  microlocalized metric $H=H_{(\la)}$.
\begin{definition}
 We call $\ric(H)$ the {\sl Ricci defect}
of the  microlocalized metric  
 $H=H_{(\la)}$.
\end{definition}

The Ricci defect of $H$ plays a fundamental role 
in the proof of the Asymptotics Theorem, see Theorem 4.5 in \cite{Einst1}  or
Theorem 2.5 in \cite{Einst2}. More precisely it appears 
as a source term in the null structure equations, see
 section 3 of \cite{Einst2}. For example the trace of 
the null second fundamental form $\chi_{AB}=H(\dd_{e_A}L,e_B)$
 $\trch=\de^{AB}\chi_{AB}$ verifies an equation, roughly, of the form
\be{ric1.11}
\frac{d}{ds}\trch=-\rr_{44}(H)+...
\end{equation}
where $\rr_{44}=\ric(L,L)=L^\a L^\b\rr_{\a\b}$ and $s$ the affine parameter
of the vectorfield $L$, i.e.  $L(s)=1$. Ignoring all other terms
on the right hand side  of \eqref{ric1.11} we see that   $\trch$
can be controlled pointwise by the mixed $L_t^1(L_x^\infty)$ norm
 of the Ricci defect. In \cite{Einst1} we have shown,
using the metric hypothesis \eqref{bootstrap} and the fact
that $H$ arises( see \eqref{ric1.3}) from an Einstein metric $\gg$,
that, 
\be{ric1.12}\|\ric(H)\|_{L_t^1L_x^\infty}\les\la^{-1-8\ep_0}.
\end{equation}

In the Asymptotics Theorem  9.1. in \cite{Einst2}  the proof of the estimates
(118-121)  was  heavily dependent on \eqref{ric1.12}.
However we also need   $L^2(\stu)$ estimates for some derivatives of 
$\trch$, in particular  the angular  derivatives $\nabb \trch$.
Differentiating the equation \eqref{ric1.11} we see that $\|\nabb\trch\|_{L^2(\stu)}$
depends on,
$$\int_u^t\|\nab \rr_{44}(H)\|_{L^2(\stau)} d\tau$$

To establish such an estimate we need  first to compare
 the Ricci defect  $\ric(H)$ with $\ric(G)=0$ and  then take advantage
of energy estimates for  derivatives of $H$ along the
null hypersurfaces $C_u$.  Here we encounter a substantial difficulty
as the $2$-surfaces  $S_{t,u}$  as well 
as the null hypersurfaces $C_u$ have been constructed  relative to the approximate
metric  $H$.  This leads to significant differences\footnote{ The estimates for the second derivatives of the higher
 frequencies of  $G$
 do in fact diverge badly.}  between the $C_u$- 
energy estimates for
 derivatives of $H$ and the corresponding ones for $G$, 
see  proposition 7.7 in
\cite{Einst2} and 
 proposition \ref{CuH} here. 

In this paper we use  the specific structure  of the component $\rr_{44}$
relative to the wave coordinates and  overcome this difficulty.
We prove the following:

\begin{theorem} On  any null hypersurface $C_u$,
\be{ricci113}
\int_u^t\|\nab \rr_{44}(H)\|_{L^2(\stau)} d\tau\les \la^{-1}
\end{equation}
\label{Ricci}
\end{theorem}
This result, stated  without proof in theorem 8.1  \cite{Einst2},  played an
essential role in the proof of the asymptotics theorem. The asymptotics theorem
itself is  a crucial step in the proof of our main theorem, see 
 \cite{Einst1}. The main goal of this paper is to prove theorem \ref{Ricci}.

\section{Preliminaries}
\subsection{Background estimates}
We start by writing down estimates for the rescaled
metric $G(t,x)=\gg(\la^{-1}t,\la^{-1}x)$. These are immediate
consequences of  
the  metric hypothesis \eqref{bootstrap} and the choice of the restricted
time interval $[0, t_*]$.
\bea
&&\|\pr G\|_{L^2_t L_x^\infty}\les
\la^{-\frac{1}2-4\ep_0},\label{asG1}\\
&&\|\pr^2 G\|_{L_t^\infty L_x^2}\les \la^{-\f12-4\ep_0}\label{asG2}
\eea
It is also easy to derive   the following estimate for $G$ in $L^2(\stu)$ norm.
\be{asG3}
\|\pr G\|_{L^2(\stu)}\les \la^{-4\eps_0}
\end{equation}
This estimate  follows by virtue of  H\"older and the trace inequality
(see theorem ?? in \cite{Einst2}) on $\stu$ from \eqref{asG2}.

We also recall the estimates for $H$ derived in
\cite{Einst1}  and \cite{Einst2}. They are
summarized in   section 7 of \cite{Einst2}.
We list below only the ones which we need
in this paper. Morally, since $H=P_{<1} G$ 
they follow from the corresponding estimates for $G$.
\begin{align}
&\|\pr H \|_{L^2_t L_x^\infty}\les
\la^{-\frac{1}2-4\ep_0},\label{aso1}\\
&\|\pr^2 H\|_{L_t^\infty L_x^2}\les \la^{-\f12-4\ep_0}
\\
&\|\pr H\|_{L^2(\stu)}\les \la^{-4\eps_0}\quad
 \label{aso2}\\
&\|\ric(H)\|_{L^1_{t} L_x^\infty}\les 
\la^{-1-8\ep_0},
\label{aso3}
\end{align} 

We also have the following cone estimates(see section 7 of \cite{Einst2} ), which play
an essential role in the proof of theorem \ref{Ricci}:
\begin{proposition}
The following estimates hold in the region
 $\Omega_*={\cal I}^+_{-1}\cap [0,t_{*}]\times
\rrrr^3$ and  $1\les\la, \mu$.
\bea
&&\|D_*\pr H\|_{L^2(C_u)} \les \la^{-\f12},\qquad
\|D_* H\|_{L^2(C_u)} \les \la^{\f12}\label{cuh'}\\
&&\|D_*\pr (P_\mu G)\|_{L^2(C_u)} \les \mu^{\frac 12 - 4\eps_0} \la^{-\frac 12 -
4\eps_0},\nn\\ &&\|D_* (P_{\mu}G)\|_{L^2(C_u)}\les \la^{-\frac 12-4\eps_0}\mu^{-\frac
12-4\eps_0}\label{10cuh'}
\eea
\label{CuH}
\end{proposition}

We shall also need  estimates for the derivatives
of the null vectorfield $L$ in $\Om_*$,
\be{nabL}
|\nab L|\les \Theta+r^{-1}
\end{equation}
where $r=r(t,u)$ is defined by
Area($S_{t,u}$)$=4\pi r^2$ and $\Theta$ verifies 
the following estimates,
\bea
\|\Theta\|_{L_t^2L_x^\infty}&\les& \la^{-\frac 12-4\ep_0}\label{ric2.1000}\\
 \|\Theta\|_{L^2(\stu)}&\les& \la^{-2\ep_0}\label{ric2.1001}
\eea
By the comparison arguments
proved in section 6.4 of \cite{Einst2} we have
\be{ric2.990}r\approx t-u
\end{equation}.
We also have,
\be{ric2.1002}
 \|\Theta\|_{L^2(\dtu)}\les \la^{-2\ep_0}
\end{equation}
where,
\begin{definition}
The annulus $\dtu$ is defined by 
$\dtu=\cup_{u\le u' \le u+1} S_{t, u'}$ is the annulus on $\Si_t$ 
of thickness $1$ and outer boundary $S_{t,u}$.
\label{defdtu}
\end{definition}
Observe that,
\be{DL}
\|\nab L\|_{L^2(\stu)}\les 1
\end{equation}
Clearly we also  have,
$$\|\nab L\|_{L^2(\dtu)}\les 1.$$

For a proof of the estimates \eqref{nabL}--\eqref{ric2.1001}
we refer to section 9 of \cite{Einst2}.

\subsection{Set-up and error terms}
\begin{definition}
We denote by $P$ the projection on the frequencies 
of size $<1$ and by $\bP$ the projection on the frequencies 
of size $\le 2 $ such  that $\bP P=P$. 
\end{definition}
\begin{definition} 
We define
\begin{eqnarray}
H(t,x)&=& P G(t,x) \nn \\
h(t,x)&=&\sum_{\mu >1} P_\mu G(t,x)\nonumber
\end{eqnarray}
\end{definition}
 Clearly,
\be{Hplush}
G=H+h
\end{equation}
Also, for the inverse metric,
\be{inverseG}
G^{-1}=(H+h)^{-1}=(I+H^{-1}h)^{-1} H^{-1}=H^{-1}-H^{-1}hH^{-1}+O(h^2)
\end{equation}
Therefore,
\bea
G_{\a\b}&=&H_{\a\b}+h_{\a\b}\label{Hplush2}\\
G^{\a\b}&=&H^{\a\b}-h^{\a\b} + O(h^2)  \label{inverseG2}
\eea
where the indices of $h$ are raised according to the matrix $H$.

In view of the fact that  $\rr_{\mu\nu}(G)=0$ we infer that,
$$
\rr_{\mu\nu}(H)=\rr_{\mu\nu}(H)-P\,\rr_{\mu\nu}(G).
$$
This is the starting point of our lengthy calculations 
 which are presented in   the following  sections. 
 In the process we 
are going to generate a large number of error terms. To better keep track
 of them we will
systematize them in the following subsection.
\subsection{Error terms}

\medn

We  start with some basic commutator estimates
which we shall need below. 
\begin{lemma} 
Let $Q$ be one of the Littlewood-Paley projections 
$Q=P, \bP, P_{\mu}$ with $\mu>1$. 
We may assume(see remark below) that the support of the integral 
kernel $Q(x)$  of the projection $Q$ is localized to 
the unit ball centered at the origin in the case $Q=P, \bP$,
and the ball of radius $\mu^{-1}$ if $Q=P_\mu$.

We denote 
$$
|Q|=\sup_{x:\,Q(x)\ne 0} |x|^{-1},
$$
Then for all $p\in [1,\infty]$ and arbitrary functions $u,w,v$
such that $\nab w, \nab f\in L^\infty_x$ and $v\in L^p_x$,
\bea
\|[Q,w]v \|_{L^p(\dtu)}&\les & \frac 
1{|Q|}\|\nab w\|_{L^\infty(\dtu)} \|v\|_{L^p(\dtu)},\label{Quv}\\
\|[Q,w]v \|_{L^p(\dtu)}&\les & \frac 1{|Q|}\|\nab w\|_{L^p(\dtu)} \|v\|_{L^\infty(\dtu)},
\label{Quv'}\\
\|[Q,\nab w]v \|_{L^p(\dtu)}&\les& \|\nab w\|_{L^\infty_x} \|v\|_{L^p_x},\label{Qduv}\\
\|\big[[Q,w],f\big ] v\|_{L^p(\dtu)}&\les& 
\frac 1{|Q|^2}\|\nab w\|_{L^\infty(\dtu)}\|\nab f\|_{L^\infty(\dtu)} 
\|v\|_{L^p(\dtu)}.\label{Q2uv}
\eea
\label{Commutation}
\end{lemma}

\begin{remark}The assumptions made on the supports of the integral
 kernels of $Q=P,\bP, P_\mu$ are essentially true\footnote{strictly
speaking they are incompatible with the compact support assumption
of the Littlewood-Paley projections in Fourier space. }. Consistent
with the uncertainty principle we can show that the kernels of 
$Q$ are rapidly decaying outside the ball of radius one for $P,\bP$
 and $\mu^{-1}$ for $P_\mu$. 
\end{remark}

\begin{proof}
The proof of the lemma is standard. For completeness  we show below  how 
 to derive estimates \eqref{Quv} and \eqref{Quv'}.
We have
\bea
[Q, w] v  &=& \il_{\Si_t} Q(x-y)\big(w(y)-w(x)\big)v(y)\,dy\nn \\ &= &
-  \il_0^1 \il_{\Si_t} Q(x-y)(x-y)^i \pr_i w(\tau y + (1-\tau)x)  
v(y)\,dy \,d\tau \label{repr}
\eea
Therefore, since the support of $Q(x)$ belongs to the unit ball 
centered at the origin,
\begin{equation}
\|[Q,w] v\|_{L^p(\dtu)}\les \frac 1{|Q|}\|\nab w\|_{L^\infty(\dtu)}
\|v\|_{L^p(\dtu)},
\label{Qwv}
\end{equation} 
where the annuli $\dtu$ on the right hand side of \eqref{Qwv}
are perhaps twice as large as the original annulus. This proves
\eqref{Quv}. 

To obtain \eqref{Quv'} we proceed as follows. Using \eqref{repr}
we obtain
\beaa
\|[Q,w] v\|_{L^p(\dtu)}&\les& \frac 1{|Q|}
\| \il_0^1 \il_{\Si_t} |Q(z)|\,\,|\nab w(x + \tau z)|\,\,  
|v(x-z)|\,dz \,d\tau \|_{L^p(\dtau)}\\&\les & 
\frac 1{|Q|}\il_0^1 \il_{\Si_t} |Q(z)| \|\nab w(\cdot\, + \tau z) v(\cdot\,
-z)\|_{L^p(\dtau)} dz\,d\tau \\ &\les &   \frac 1{|Q|} \|\nab w\|_{L^p(\dtau)}
\|v\|_{L^\infty (\dtau)}.
\eeaa
as desired.
Here we once again used that  the support of $Q(z)$ belongs to the unit ball 
centered at the origin.
\end{proof}

\begin{definition} 
Given functions $f,v,w$ in $ L^{\infty}(\Om_*) $
we introduce the following:

\begin{itemize}
\item 
We denote by $[f]$ any operator with the 
property that
for any function $v$ in $ \Om_*$ and any $t\in[0,t_*]$, $u\ge -1$:
\be{pione}
\|[f]\cdot v\|_{L^2(\dtu)}\les
\min\bigg\{\|f\|_{L^\infty(\dtu)}
\|v\|_{L^2(\dtu)}\,,\, \|f\|_{L^2(\dtu)}
\|v\|_{L^\infty(\dtu)}\bigg\}
\end{equation}
\item We denote by $\pi(f, v; w)$ any function in $ \Om_* $ which
satisfies the inequality:
\be{pithree}
\|\pi(f,\,v\,;\,w)\|_{L^2(\dtu)} \les
\|f\|_{L_x^\infty(\dtu)}\,\|v\|_{L_x^\infty(\dtu)} \|w\|_{L^2(\dtu)}
\end{equation}
\end{itemize}
\label{ric2.989}
\end{definition}
\begin{definition}
Given two operators $A$ and $B$ we say that $A\les B$  if for any function $v$
\be{alesb}
\|A v\|_{L^2(\dtu)}\les \|B v\|_{L^2(\dtu)}
\end{equation}
We also say that $\pi(f,v;w)\les \pi(g,v;w)$ if 
\be{pfg}
\|\pi(f,\,v\,;\,w)\|_{L^2(\dtu)} \les
\|g\|_{L_x^\infty(\dtu)}\,\|v\|_{L_x^\infty(\dtu)} \|w\|_{L^2(\dtu)}
\end{equation}
\label{defles}
\end{definition}

\begin{remark} The expression $[f]$
verifies the following  trivial property $$[af]\les \|a\|_{L^\infty}[f].$$
The same holds true for   $\pi(f,g;v)$ with respect to all entries.
\end{remark}
For the Littlewood-Paley projection $Q$ let $(Qf)$ be the result
of the application of $Q$ to $f$. We also denote by $Q f$ the operator 
whose action on functions is defined by
$$
Qf (v): = Q(fv)
$$
The typical examples of 
expressions of type $[f]$ are listed in the following lemma.
\begin{lemma} 

\medn
\begin{itemize}
\item  For the projections $Q=I, P, \bP, P_\mu$  with dyadic $\mu >1$,  
we have $Q f\les [f]$ and
$(Q f)\les [f]$.
\item For $Q=P, \bP, P_\mu$, we have, 
$[Q,f]\les \frac 1{|Q|}[\nab f]$ and 
$[\nab Q, f] \les [\nab f]$. 
\end{itemize}
\label{propbr}
\end{lemma}
\begin{proof}
To verify that $Q f\les [f]$ we estimate
\beaa
\|Q f (v)\|_{L^2(\dtu)} &\les& \|Q (f v)\|_{L^2(\dtu)}\les
\|f v\|_{L^2(\dtu)}\\ &\les&  \min\bigg\{\|f\|_{L^\infty(\dtu)}
\|v\|_{L^2(\dtu)}\,,\, \|f\|_{L^2(\dtu)}
\|v\|_{L^\infty(\dtu)}\bigg\}
\eeaa
A similar argument shows that $(Qf)\les [f]$.
We now verify that $[Q,f]\les \frac 1{|Q|}[\nab f]$.
Using the commutator estimates \eqref{Quv} and \eqref{Quv'}
we obtain 
$$
\|[Q,f] v\|_{L^2(\dtu)} \les
\frac 1{|Q|}\min\bigg\{\|\nab f\|_{L^\infty(\dtu)}
\|v\|_{L^2(\dtu)}\,,\, \|\nab f\|_{L^2(\dtu)}
\|v\|_{L^\infty(\dtu)}\bigg\}
$$
as desired.
The proof of the estimate $[\nab Q, f] \les [\nab f]$ is similar.
It uses the commutator estimate \eqref{Qduv}.
\end{proof}

We also record some similar  properties of the triple expressions $\pi(f,g;h)$.
\begin{lemma}

\medn
\begin{itemize}
\item For $Q_i=I, P, \bP, P_\mu$ with some dyadic $\mu\ge 1$ and $i=1,..,3$, we have 
\begin{equation}
(Q_1 f) (Q_2 v) (Q_3 w) \les \pi(f,v ; w).
\label{tric3.1}
\end{equation}
\item With the same choice of $ Q_1, Q_2$,
\begin{align}
[f] \bigg((Q_1 v) (Q_2 w)\bigg) &\les  \pi (v, w; f),\label{tripleform}\\
[f] \bigg((Q_1 v) (Q_2 w)\bigg) &\les  \pi (f, v; w).\label{tripleform'}
\end{align}
\item With the same choice of $ Q$
\begin{align}
[f] [v] (Q w) &\les  \pi (f, v; w),\label{triform}\\
[f] [v] (Q w) &\les  \pi (v, w; f).\label{triform'}
\end{align}
\item If $\|f\|_{L^\infty}\les \|g\|_{L^\infty}$  then
\begin{equation}
\pi (f,v; w)\les \pi(g,v;w)
\label{tric3.4}
\end{equation}
\end{itemize} 
\label{triplepi}
\end{lemma}
\begin{proof}
The proof of \eqref{tric3.1} follows immediately from the definition 
of $\pi(f,v;w)$ and the properties of the projection $Q_i$. Indeed
\beaa
\|(Q_1 f) (Q_2 v) (Q_3 w)\|_{L^2(\dtu)} &\les& 
\|(Q_1 f)\|_{L^\infty(\dtu)}\|(Q_2 v)\|_{L^\infty(\dtu)}
\|(Q_3 w)\|_{L^2(\dtu)}\\ &\les& 
\|f\|_{L^\infty(\dtu)}\|v\|_{L^\infty(\dtu)}
\| w\|_{L^2(\dtu)}
\eeaa
To obtain \eqref{tripleform} we estimate using 
definition \eqref{pione} for $[\,\,]$,
\beaa
\|[f] \bigg((Q_1 v) (Q_2 w)\bigg)\|_{L^2(\dtu)}&\les &
\|f\|_{L^\infty(\dtu)} \|(Q_1 v) (Q_2 w) \|_{L^2(\dtu)}\\
&\les & \|f\|_{L^\infty(\dtu)}\|v\|_{L^\infty(\dtu)} 
\|w\|_{L^2(\dtu)} 
\eeaa
The alternative estimate in \eqref{pione} for $[f]$ similarly leads
to \eqref{tripleform'}.

We also derive 
\beaa
\|[f] [v] (Q w)\|_{L^2(\dtu)} &\les & 
\|f\|_{L^\infty(\dtu)}\| [v] (Q w)\|_{L^2(\dtu)}\\
 &\les & \|f\|_{L^\infty(\dtu)} \|v\|_{L^\infty(\dtu)}
\| w)\|_{L^2(\dtu)}
\eeaa
as claimed in \eqref{triform}. The estimate \eqref{triform'} 
once again is obtained by using the alternative term in \eqref{pione} 
in the estimates for $[f]$ and $[v]$.  

Finally, if 
$\|f\|_{L^\infty(\dtu)}\les \|g\|_{L^\infty(\dtu)}$,
then,
\beaa
\|\pi(f,v;w)\|_{L^2(\dtu)} &\les & 
\|f\|_{L^\infty(\dtu)} \| v \|_{L^\infty(\dtu)}
\| w\|_{L^2(\dtu)}\\
&\les & \|g\|_{L^\infty(\dtu)} \| v \|_{L^\infty(\dtu)}
\| w\|_{L^2(\dtu)}
\eeaa
Thus according to definition \ref{defles}, 
$\pi(f,v;w)\les \pi(g,v;w)$ as desired in \eqref{tric3.4}.
\end{proof}

\section{Wave coordinate condition}

In what follows we shall  rely crucially  on the fact that
our standard coordinates $x^\a,\,\a=0,..,3$ satisfy the wave
coordinate condition \eqref{ric1.2} relative to the metric $G$.
Recall that the wave coordinate condition has the form:
\beaa
0&=&\frac{1}{\sqrt{|G|}}\pr_\a(G^{\a\b}\sqrt{|G|})\\
&=&\pr_\a G^{\a\b}+\f12  G^{\a\b}G^{\ga\de}\pr_\a G_{\ga\de}
\eeaa
or, in view of $\pr(G^{\a\b}G_{\b\si})=0$,
\be{Gwave}
G^{\a\b} \pr_\a G_{\b\si}=\half G^{\a\b}\pr_\si G_{\a\b}.
\end{equation}

Next we shall review some basic notation
connected to our standard null frame $L=e_{4}$, $\Lb=e_3$, $e_A$, $A=1,2$.
When $L,\Lb$ are applied to scalar quantities we also use the notation
$L=\partial_4$, $\Lb=\pr_3$.
Recall that the null components of the metric $H$ are given by,
$$H_{34}=-2, \quad H_{33}=H_{44}=H_{3A}=H_{4A}=0,\quad H_{AB}=\de_{AB}.$$
The  null components of the  inverse metric are therefore,
$$H^{34}=-\f12, \quad H^{33}=H^{44}=H^{3A}=H^{4A}=0,\quad H^{AB}=\de^{AB}.$$
Given a vectorfield $X=X^\a\pr_\a$ we decompose relative to the null frame as
follows:
\bea
X&=&-\f12 <L,X>\Lb -\f12 <\Lb,X> L +<e_A, X>e_A\nn\\
&=&-\f12 X_4\Lb -\f12 X_3 L +X_Ae_A \label{Xdown}
\eea
or, using upper indices,
\be{Xup}
X=X^3\Lb+X^4L +X^A e_A
\end{equation}
where
$$X^3=-\f12 X_4,\quad X^4=-\f12 X_3,\quad X^A=X_A.$$
In view of this we shall use the following notation,
\begin{definition}
\label{nullframePi}

For an arbitrary  spacetime tensor $M^{\a\b}$,
$$
\aligned
\ &M^{3\b}:=-\frac 12 M^{\a\b} L_\a = 
 -\frac 12 M^{\a\b} H_{\a\ga}L^\ga=-\f12 M^{\a\b}H_{\a4} ,\\
\ &M^{4\b}:=-\frac 12 M^{\a\b} \Lb_\a =  
-\frac 12 M^{\a\b} H_{\a\ga}\Lb^\ga= -\f12M^{\a\b}H_{\a3},\\
\ &M^{A\b}:=M^{\a\b} e_{A\a} =  M^{\a\b} H_{\a\ga}e_{A}^\ga
\endaligned
$$ 
In particular
$$H^{3\a}=H^{34}L^\a=-\f12 L^\a. $$
\end{definition}

\begin{definition}
Given a scalar function  $f$ we shall denote by $\DC f$ any
function for which we have an estimate of the
form
$$|\DC f|\les
|L(f)|+(\,\sum_{A=1,2}|e_A(f)|^2\,)^{\f12}=|L(f)|+|\nabb f|$$
Given a tensorfield  $U$ with components
 $U_{\underline{\a}}^{\,\,\,\underline{\b}}$ relative to our standard
coordinates $x^\a$ we denote by $\DC U$ a scalar quantity 
which can be estimated by,
$$|\DC U|\les \sum_{\underline{\a},\underline{\b}}|\DC
U_{\underline{\a}}^{\,\,\,\underline{\b}}|. $$
Given two tensors $U, V$  we denote by $U\DC V$  a scalar  quantity 
which can be estimated by 
$$|U\DC V|\les |U||\DC V|.$$
\label{tangentialnotation}
\end{definition}

For example, consider  the coordinate vectorfield $\pr_\a$ and 
decompose it relative to the null frame $L,\Lb,e_A$ according to
\eqref{Xdown}. We shall write the decomposition formula in the form,
\be{decompartial}
\pr_\a=-\f12 L_\a \, \Lb +\DC. 
\end{equation}

Using the above notation we are now ready to state 
the main result of this section.
\begin{lemma}
The following identities\footnote{They are
in fact approximate identities. The terms on the right hand
side are schematically. What we mean is that the terms
on the left can be estimated by the quantities
appearing on the right.} are consequencies of the wave
coordinate condition \eqref{Gwave}:
\bea
2 H^{3\a} \pr_3 (QG)_{\a\si}&=&  H^{\a\b} \pr_\si (QG)_{\a\b} +G\DC(Q G) +\err
\label{2mod11}\\
L^\mu L^\nu \partial_3(Q G_{\mu\nu}) &=& G\cdot \DC(Q G) + \err,\label{2mod1000}\\
L^\a e_A^\si \pr_3(  QG_{\a\si}) &=& G\cdot\DC (QG) + \err,
\label{2mod13}\\
\err&=&h\pr (Q G) +\frac{1}{|Q|}[\pr G]\pr G
\eea
In particular,
\be{2mod12}
L^\a L^\si \pr_3  H_{\a\si} = G\cdot \DC H + \err
\end{equation}
We also have,
\bea
L^\si L^\b \pr_\a \pr_\ga (QG_{\b\si}) &= &
H \cdot \DC \pr (QG) + \err\label{2mod14}\\
\err&=&\frac{1}{|Q|}\bigg([\pr G]\c\pr^2 G+[\pr^2 G]\pr G\bigg)
+h\c\pr^2 (QG) 
+ \pr G\c\pr(QG)\nn
\eea

\label{Wcord}
\end{lemma}
\begin{proof}
We start by projecting \eqref{Gwave}:
$$ Q\big(G^{\a\b} \pr_\a G_{\b\si}\big)=\half Q\big( G^{\a\b}\pr_\si G_{\a\b}\big).$$
In view of the fact that $Q(u\cdot v)=u\cdot Q v +\frac{1}{|Q|}[\nab u]v$
we derive
\be{ric3.1}
G^{\a\b} \pr_\a (QG_{\b\si})=\half G^{\a\b}\pr_\si (QG_{\a\b})+
\frac{1}{|Q|}[\pr G]\pr G
\end{equation}
Expanding \eqref{ric3.1} relative to the null frame, we have 
$$ G^{3\b}\pr_3 (QG_{\b\si})+ G^{4\b}\pr_4(QG_{\b\si}) +
G^{A\b}\pr_A(QG_{\b\si})=\half
G^{\ga\de}\pr_\si (QG_{\ga\de})+\frac{1}{|Q|}[\pr G]\pr G,$$
whence, for any $\si$, 
$$ G^{3\b}\pr_3(QG_{\b\si})=\half
G^{\ga\de}\pr_\si (QG_{\ga\de})+G\DC( QG)+\frac{1}{|Q|}[\pr G]\pr G
$$

Writing $G^{3\b} = H^{3\b} - h^{3\b} + O(h^2)$ 
we derive, 
\bea
2 H^{3\a} \pr_3 (QG_{\a\si}) &=& H^{\a\b} \pr_\si(Q G_{\a\b}) +G\DC( QG) +\err
\label{ric3.2}\\
\err&=& h\pr(Q G) +h^2\pr(Q G)+\frac{1}{|Q|}[\pr G]\pr G
\eea
\begin{remark}
 Since $\|h\|_{L^\infty}\les 1$, the error term $h^2\pr(Q G)$
can be treated in the same way as $h\pr(Q G)$ and we shall ignore it.
In what follows we shall often drop
terms like this without  further mentioning.
\end{remark}
 We thus derive the desired 
approximate identity \eqref{2mod11}. 

Contracting  \eqref{ric3.2} with $L^\si$ 
we obtain,
$$
 2L^\si\,H^{3\a}\pr_3(QG_{\a\si})=G\DC(Q G) + H\DC(Q G)+ \err
$$
As $H\DC G$ can be estimated exactly  in the same way as the more difficult
term $G\DC G$ we shall drop it. We shall later absorb similar
terms into related, more difficult terms, without further
mentioning.

We now recall that 
$H^{3\a}=-\f12 L^\a$.
Henceforth,
$$
- L^\nu L^\mu \,\pr_3(QG_{\mu\nu})= G\DC( QG) + \err,
$$
which gives \eqref{2mod1000}.

We can also contract \eqref{ric3.2} with $e_A^\si$ to obtain
$$
e_A^\si\,H^{3\a}\pr_3(QG_{\a\si})=G\DC (QG) + \err.
$$ 
 Using again the relation
$H^{3\a}=-\f12 L^\a$,   \eqref{2mod13} immediately follows. 
 We shall now prove \eqref{2mod14}. Differentiating \eqref{Gwave},
 we find,
\be{8.101}
G^{\a\b} \pr_\ga\pr_\a G_{\b\si}=\half G^{\a\b}\pr_\ga \pr_\si G_{\a\b} 
+ \pr G\c\pr G.
\end{equation}
We manipulate the left
hand side of
\eqref{8.101} schematically as follows:
\beaa
Q(G\cdot \pr^2 G)&=&H\c\pr^2 (QG)+Q(G\cdot \pr^2 G)-H\c\pr^2 (QG)\qquad\qquad\qquad
\,\,\,\,\,\\ &=&H\c\pr^2 (QG)+Q(G\cdot \pr^2 G)-G\c\pr^2( QG) +h\c\pr^2 (QG)
\qquad (\mbox{as}\,\, G=H-h)\\
&=&H\c\pr^2 (QG)+[Q,G]\pr^2 G+h\c\pr^2 (QG)\\
&=&H\c\pr^2 (QG)+\frac{1}{|Q|}[\pr G]\pr^2 G+h\c\pr^2( QG)
\eeaa
Therefore, proceeding in the same way on the right hand
side of \eqref{8.101},
\bea
H^{\a\b} \pr_\ga\pr_\a (QG_{\b\si}) &=& \half Q\bigg(G^{\a\b}\pr_\ga \pr_\si
G_{\a\b}\bigg) +\err\nn\\
&=& \half H^{\a\b}\pr_\ga \pr_\si
(QG_{\a\b}) +\err\label{Lsi}\\
\err&=&\frac{1}{|Q|}[\pr G]\c\pr^2 G+h\c\pr^2 (QG)+ Q(\pr G\c\pr G)\nn\\
&=&\frac{1}{|Q|}[\pr G]\c\pr^2 G+h\c\pr^2 (QG)+\frac{1}{|Q|}[\pr^2 G]\pr G 
+ \pr G\c\pr(QG)\nn
\eea

We now  contract \eqref{Lsi} with  $L^\si$.

\beaa
L^\si H^{\a\b} \pr_\ga\pr_\a (QG_{\b\si})& = &\half L^\si
  H^{\a\b}\pr_\si\pr_\ga 
(QG_{\a\b}) + 
\err\\
& = &H\c \DC\pr ( QG )  +\err
\eeaa
Therefore, expressing $H^{\a\b}\pr_\a$ relative to the null
frame $L,\Lb,e_A$, 
$$
L^\si L^\b \pr_\a \pr_\ga (QG_{\b\si}) = 
H \cdot \DC \pr (QG) + \err
$$
\end{proof}
\section{First reduction} 
In this section we show how to reduce the statement of Theorem \ref{Ricci}
to the following:
\be{R0}
\int_{u+1}^t \|\rr_{44}(H)\|_{L^2({\dtau})}d\tau\les \la^{-1}
\end{equation}
where, see definition \ref{defdtu},   $\dtau=\cup_{u\le u' \le u+1} S_{\tau, u'}$ is the
annulus on
$\Si_\tau$  of thickness $1$ and outer boundary $S_{\tau,u}$. Throughout
this  and the remaining sections  we denote $\rr_{44}=\rr_{44}(H)$ and $\ric=\ric(H)$.

\medn
{\em Step 1}\quad  Take care of $\int_u^{u+1}\|\nab\rr_{44}\|_{L^2({\stau})}d\tau$.
\medn

 We  start with formula 
\be{R1}
\nab \rr_{44}=\nab L^\mu L^\nu \rr_{\mu\nu}=2(\nab L) \cdot L\cdot \ric+  L^\mu L^\nu
\nab \rr_{\mu\nu}.
\end{equation}
Recall  that, see \eqref{nabL},
$$\nabla L\les r^{-1}+\Theta$$ with $\Theta$ verifying the estimates 
\eqref{ric2.1000}--\eqref{ric2.1001}. Clearly, 
 $\|(\Theta + r^{-1})\|_{L^2({\stau})}\les 1$.
Also observe that, in view of \eqref{ric2.990},  $r\le 1$ as $\tau$ varies between $u$
and
$u+1$.
 We infer that
$$
\aligned
\int_u^{u+1}\|\nab\rr_{44}\|_{L^2({\stau})}d\tau &\les 
\int_u^{u+1} \bigg(\|(\Theta + r^{-1})\|_{L^2({\stau})}\|\ric \|_{L^\infty_x} +
r \|\nab \ric \|_{L^\infty_x}\bigg) d\tau \\ &\les \|\ric\|_{L^1_t L^\infty_x}
 + \|\nab \ric\|_{L^1_t L^\infty_x}.
\endaligned
$$
 It remains to observe that 
the frequencies of $\ric(H)$ are essentially $\le 2$ and therefore,
$\|\nab\ric\|_{L^\infty_x}\les \|\ric \|_{L^\infty_x}$.
Henceforth, in view of the background estimate \eqref{aso3},
\be{R2'}
\int_u^{u+1}\|\nab\rr_{44}\|_{L^2({\stau})}d\tau \les 
 \|\ric\|_{L^1_t L^\infty_x}\les \la^{-1-4\ep_0}
\end{equation}

{\em Step 2}\quad Take care of $\int_{u+1}^t\|\nabla\rr_{44}\|_{L^2({\stau})}d\tau$.

We start as in Step 1 with  formula \eqref{R1}.
To  estimate the first term on the right  hand side  of \eqref{R1} we
use 
 $\|\nabla L\|_{L^2(\stau)}\les 1$, see \eqref{DL}.

Therefore, using also \eqref{aso3}
\beaa\int_{u+1}^t\|\nabla L\cdot L\ric\|_{L^2({\stau})}d\tau&\les& 
\sup_{u+1\le \tau\le t}\|\nabla L\|_{L^2(\stau)}\|\ric\|_{L_t^1L_x^\infty}
\les \la^{-1-4\ep_0}
\eeaa
as desired. In other words,
\be{R1'}
\int_{u+1}^t\|\nabla\rr_{44}\|_{L^2({\stau})}d\tau\les \la^{-1-4\ep_0}+
\int_{u+1}^t\| L^\mu L^\nu \nab\rr_{\mu\nu}\|_{L^2({\stau})}d\tau
\end{equation}

 It remains to estimate the second term in \eqref{R1'}.
Using the simple estimate:
$$\|f\|_{L^2(\stau)}^2\les \|\nabla f\|_{L^2(\dtau)}\|f\|_{L^2(\dtau)}$$
where $\dtau=\cup_{u\le u' \le u+1} S_{\tau, u'}$ is the annulus on $\Si_\tau$ 
of thickness $1$ and outer boundary $S_{\tau,u}$.
\beaa
\| L^\mu L^\nu \nab\rr_{\mu\nu}\|_{L^2({\stau})}&\les& \|\nabla L^\mu L^\nu
\nab\rr_{\mu\nu}\|_{L^2({\dtau})}^{\f12}\| L^\mu L^\nu
\nab\rr_{\mu\nu}\|_{L^2({\dtau})}^{\f12}\\ &\les&  \|\nabla L^\mu L^\nu
\nab\rr_{\mu\nu}\|_{L^2({\dtau})}+\| L^\mu L^\nu
\nab\rr_{\mu\nu}\|_{L^2({\dtau})}
\eeaa

Now, using $\|\nab L\|_{L^2({\dtau})}\les 1$,
\bea
\|\nabla L^\mu L^\nu
\nab\rr_{\mu\nu}\|_{L^2({\dtau})}&\les& \|\nab L\c\nab \ric\|_{L^2({\dtau})} +
\|L^\mu L^\nu
\nab^2\rr_{\mu\nu}\|_{L^2({\dtau})}\nn\\ &\les& 
\|\nab L\|_{L^2({\dtau})}\|\nab \ric\|_{L^\infty_x} + \|L^\mu L^\nu
\nabla^2 \rr_{\mu\nu}\|_{L^2({\dtau})}\nn\\
&\les& \|\nab \ric\|_{L^\infty_x} + \|L^\mu L^\nu
\nabla^2 \rr_{\mu\nu}\|_{L^2({\dtau})}
\label{R4}
\eea
Now, since $\rr_{\mu\nu}\approx \bP \rr_{\mu\nu}(H)$,
and $\|\nabla^m \bP f\|_{L^2({\dtau})}\les \|f\|_{L^2({\dtau})}$
with  perhaps a  slightly larger annulus $\dtau$,

\beaa
\|L^\mu L^\nu
\nabla^2 \rr_{\mu\nu}\|_{L^2({\dtau})}&\les& \|\nabla^2\bP L^\mu L^\nu
 \rr_{\mu\nu}\|_{L^2({\dtau})}+\|\big[L^\mu L^\nu,
\nabla^2\bP\big]\ric\|_{L^2({\dtau})} \\
&\les&\| L^\mu L^\nu
 \rr_{\mu\nu}\|_{L^2({\dtau})}+\|\big[L^\mu L^\nu,
\nabla^2\bP\big]\ric\|_{L^2({\dtau})}
\eeaa
To treat the second term we shall use the following 
commutation lemma,\, see lemma \ref{Commutation},
\be{comm1}
\|\big[ f, \nabla^k \bP\big] g\|_{L^2(\dtau)}
\les\|\nabla f\|_{L^2(\dtau)}\|g\|_{L^\infty_x}
\end{equation}
with a possible larger annulus $\dtau$ on the right hane side.

Therefore,
\beaa\|\big[L^\mu L^\nu,
\nabla^2\bP\big]\ric\|_{L^2({\dtau})}&\les &\|\nabla L\|_{L^2(\dtau)}
\|\ric\|_{L_x^\infty}\les \|\ric\|_{L_x^\infty}
\eeaa
Therefore, back to \eqref{R4},
$$
\|\nabla L^\mu L^\nu
\nabla \rr_{\mu\nu}\|_{L^2({\dtau})}\les \|
 \rr_{44}\|_{L^2({\dtau})}+\|\ric\|_{L_x^\infty}+\|\nabla \ric\|_{L_x^\infty}
$$
or, since $\ric=\ric(H)\approx \bP \ric $,
\be{R6}
\|\nabla L^\mu L^\nu
\nabla \rr_{\mu\nu}\|_{L^2({\dtau})}\les \|
 \rr_{44}\|_{L^2({\dtau})}+\|\ric\|_{L_x^\infty}
\end{equation}
Also, clearly,
$$
\| L^\mu L^\nu
\nabla \rr_{\mu\nu}\|_{L^2({\dtau})}\les \|
 \rr_{44}\|_{L^2({\dtau})}+\|\ric\|_{L_x^\infty}.
$$
Therefore,
\be{R7}
\| L^\mu L^\nu \nab\rr_{\mu\nu}\|_{L^2({\stau})}
\les \|\rr_{44}\|_{L^2({\dtau})}+\|\ric\|_{L_x^\infty}
\end{equation}
whence,
\beaa
\int_{u+1}^t\|\nab\rr_{44}\|_{L^2({\stau})}d\tau
&\les&\int_{u+1}^t\|\rr_{44}\|_{L^2({\dtau})}d\tau+
\|\ric\|_{L_t^1L_x^\infty}+\la^{-1-4\ep_0}
\\
&\les&\int_{u+1}^t\|\rr_{44}\|_{L^2({\dtau})}d\tau+\la^{-1-4\ep_0}
\eeaa 
Combining this with \eqref{R2'} we obtain,
\be{R10}
\int_u^t\|\nab \rr_{44}(H)\|_{L^2(\stau)} d\tau
\les\int_{u+1}^t \|\rr_{44}\|_{L^2({\dtau})}d\tau+\la^{-1-4\ep_0} 
\end{equation}
as desired.

\section{The algebraic structure of  $\rr_{\mu\nu}(H)$}
\label{algebraicstructure}
We start with the formula,
\be{ricci100}
\rr_{\mu\nu}(H)=\rr_{\mu\nu}(H)-P\,\rr_{\mu\nu}(G)
\end{equation}

Recall the expression of the  Ricci tensor relative to local coordinates:
\begin{eqnarray}
\rr_{\mu\nu}(H)&=&\rra_{\mu\nu}(H)+\rrb_{\mu\nu}(H)\label{ricci99}\\
\rra_{\mu\nu}(H)&=&\frac 12 H^{\a\b}\bigg(H_{\a\nu\,, \b\mu}+H_{\b\mu\,,
\a\nu}-H_{\a\b\,, \mu\nu}-H_{\mu\nu\,, \a\b}\bigg)\nn\\
&=&: \frac 12 H^{\a\b}
H_{[\a\b\mu\nu]}= \frac 12 H^{\a\b}\big(
H_{[[\a\b\mu\nu]]}-H_{\mu\nu\,,\a\b}\big) \nonumber\\
\rrb_{\mu\nu}(H)&=&H^{\a\b} H_{\ga\de}\bigg(\Ga^\ga_{\mu\b}(H)
\Ga^\de_{\a\nu}(H)
-\Ga^\ga_{\mu\nu}(H)\Ga^\de_{\a\b}(H)\bigg)\nonumber
\end{eqnarray}
where
$$\Ga^{\ga}_{\a\b}(H)=\frac{1}{2}H^{\ga\si}\bigg( H_{\si\b\,, \a}+
 H_{\a\si\,, \b}- H_{\a\b\,, \si}\bigg).
$$
To calculate $\rr_{\mu\nu}(H)- P\,\rr_{\mu\nu}(G)$ we use \eqref{Hplush2}
and  \eqref{inverseG2},
\bea
G_{\a\b}&=&H_{\a\b}+h_{\a\b}\nn\\
G^{\a\b}&=&H^{\a\b}-h^{\a\b} + [h\c h].\nn
\eea
Therefore, using the notation in \eqref{ricci99} and the fact 
that $H= P\,G$,
we find,
\be{ricci101}
\begin{split}
\rra_{\mu\nu}(H)-P\,\rra_{\mu\nu}(G)&=\frac{1}{2}\bigg(H^{\a\b}H_{[\a\b\mu\nu]}-
P \,(G^{\a\b}G_{[\a\b\mu\nu]})\bigg)\\
&=\frac{1}{2}\bigg(h^{\a\b}H_{[\a\b\mu\nu]} + G^{\a\b}H_{[\a\b\mu\nu]}
- P\,(G^{\a\b}G_{[\a\b\mu\nu]})\bigg) +
h^2 \pr^2H \nonumber \\
&=\frac{1}{2}\bigg (h^{\a\b}H_{[\a\b\mu\nu]} + G^{\a\b}P\,G_{[\a\b\mu\nu]}
- P\,(G^{\a\b}G_{[\a\b\mu\nu]})\bigg) +
h^2 \pr^2 H \nonumber\\
&=\frac{1}{2}\bigg (h^{\a\b}H_{[\a\b\mu\nu]} + [G^{\a\b},P]\,G_{[\a\b\mu\nu]}
\bigg) +  \pi(h,h;\pr^2 H)
\end{split}
\end{equation}
For convenience we shall introduce the following
notation,
\begin{definition}
Given two  scalar functions $v,w$ we define
$$\{v\,, \,w\}'=[v, P]\cdot w.$$
\end{definition}
Therefore,
\be{IandII}
\rra_{\mu\nu}(H)-P\,\rra_{\mu\nu}(G)=\frac{1}{2}\bigg (h^{\a\b}H_{[\a\b\mu\nu]} 
+ \{G^{\a\b}\,,\, G_{[\a\b\mu\nu]}\}'
\bigg) +  \pi(h,h;\pr^2 H)
\end{equation}
\begin{remark} Observe that,
\begin{equation}
\begin{split}
[P,v]( I-\bP)w &= P\big (v(I-\bP)w\big)-vP(I-\bP)w\big)\\
& = \sum_{\la_1>1,\la_2>2, |\ln (\la_1\la_2^{-1})|\le 2} 
P \big(v^{\la_1} w^{\la_2}\big)\\
&= \sum_{\la_1>1,\la_2>2, |\ln (\la_1\la_2^{-1})|\le 2} 
[P, v^{\la_1}] w^{\la_2}
\end{split}
\label{10mod13}
\end{equation}
Thus, writing $w=\bP w+(I-\bP )w $,
\bea
\{v\,, \,w\}'
&=& [v, P]\bP w +
 \sum_{\la_1>1,\la_2>2, |\ln (\la_1\la_2^{-1})|\le 2} 
[v^{\la_1}, P] w^{\la_2}
\label{PGPG} 
\eea
\end{remark}

To compute the contribution
to \eqref{ricci100} of the quadratic terms $\rr_{\mu\nu}^{(2)}(H)$
we start $\Ga^\ga_{\a\b}(H)$ which we write in the form
$$
\Ga^\ga_{\a\b}(H) = \frac{1}{2}H^{\ga\si}\bigg( P\,G_{\si\b\,,
\a}+ P\,G_{\a\si\,, \b}- P\,G_{\a\b\,, \si}\bigg). $$
Now  commuting $P$ with
$H$, and using \eqref{inverseG2},
$$
\Ga^\ga_{\a\b}(H) = 
P\,\Ga^\ga_{\a\b}(G) + [\pr H]\pr G + [h \pr G] 
$$
Therefore,  using that $h$ and $\pr H$ are bounded 
 and the definition of the error term $\pi$, we infer that 
\bea
\rrb_{\mu\nu}(H) &=& H^{\a\b} H_{\ga\de}\bigg(\big(P\,\Ga^\ga_{\mu\b}(G)\big) 
\big (P\,\Ga^\de_{\a\nu}(G)\big)
-\big(P\,\Ga^\ga_{\mu\nu}(G)\big) \big (P\,\Ga^\de_{\a\b}(G)\big) \bigg) \nn\\& +& 
\pi(\pr H,\pr H;\pr G) + \pi(h,\pr G;\pr G) 
\label{R2H}
\eea
On the other hand, using first the formulae \eqref{Hplush2}, \eqref{inverseG2}
and then commuting $P$ with $H$,
\bea
P\,\rrb_{\mu\nu}(G)&=&P\bigg( G^{\a\b} G_{\ga\de}\, \big(\,\,\Ga^\ga_{\mu\b}(G)
\Ga^\de_{\a\nu}(G) -\Ga^\ga_{\mu\nu}(G)\Ga^\de_{\a\b}(G)\,\,\big)\bigg)\nn\\
 &=&  H^{\a\b} H_{\ga\de}\,P \bigg(\Ga^\ga_{\mu\b}(G)
\Ga^\de_{\a\nu}(G) -\Ga^\ga_{\mu\nu}(G)\Ga^\de_{\a\b}(G)\bigg)\nn\\
&+& 
\pi(\pr H,\pr G, \pr G) + \pi(h,\pr G,\pr G)\label{PR2G}
\eea
Thus, combining \eqref{R2H} with \eqref{PR2G},
\bea
\rrb_{\mu\nu}(H) - P\,\rrb_{\mu\nu}(G) &=& - H^{\a\b} H_{\ga\de}
\bigg (P\, \big(\Ga^\ga_{\mu\b}(G)\Ga^\de_{\a\nu}(G)\big )- 
\big(P\,\Ga^\ga_{\mu\b}(G)\big) \big (P\,\Ga^\de_{\a\nu}(G)\big)\nn\\
&-& P\,\big(\Ga^\ga_{\mu\nu}(G)\Ga^\de_{\a\b}(G)\big) + 
\big(P\,\Ga^\ga_{\mu\nu}(G)\big) \big (P\,\Ga^\de_{\a\b}(G)\big)\bigg)
\nn\\
 &+&\pi(\pr H,\pr H,\pr G) + \pi(h,\pr G,\pr G)
\label{R2H-PR2G}
\eea
To simplify the expression above we  introduce the following,
\begin{definition}\label{paradifferential} Given  two functions $v$ and $w$
 we introduce their \textit{modified\footnote{It differs from the standard
paradifferential product. In our definition we have removed 
the low-low interactions. } paradifferential product} 
$\{v,w\}$.
\be{lowlowremoved}
\{v,w\}:=P(v\cdot w)-Pv\cdot Pw
\end{equation}
\end{definition}
\begin{remark}
Observe that,
\begin{align}
\{v,w\}&= P\,\sum_{\half<\la_1\le 4 }
v^{\la_1} P_{\le \f12}\, w
+ P\,\sum_{\la_1>\half, |\ln (\la_1\la_2^{-1})|\le 2}
v^{\la_1}\,  w^{\la_2}\nn  
 \\ &+ P\,\sum_{\half < \la_2\le 4 }
P_{\le \f12} v\, w^{\la_2} 
+ P\,\sum_{\la_2 >\half, |\ln (\la_1\la_2^{-1})|\le 2}
v^{\la_1}\, w^{\la_2}\label{figdec}\\
& - \sum_{\half<\la_1\le 1 }
v^{\la_1}\, P_{\le \f12}w
- \sum_{\half <\la_2\le 1}
P_{\le \f12}v\,w^{\la_2}\nn  
 \\ &- \sum_{\half < \la_1, \la_2\le 1 }
v^{\la_1} \, w^{\la_2}\nn
\end{align}
\label{remark5.5}
\end{remark}

With this definition  we can write
\begin{align}
\rrb_{\mu\nu}(H) - P\,\rrb_{\mu\nu}(G) &= - H^{\a\b} H_{\ga\de}
\bigg(\big\{\Ga^\ga_{\mu\b}(G), \Ga^\de_{\a\nu}(G)\big\}\nn \\ &- 
\big\{\Ga^\ga_{\mu\nu}(G),\Ga^\de_{\a\b}(G)\big\}\bigg)
+\err
\label{2mod1}
\end{align}
with the error term of the form
$$
\err = \pi(\pr H,\pr H,\pr G) + \pi(h,\pr G,\pr G)
$$

Thus,  taking into account \eqref{IandII} and 
\eqref{2mod1}, we  rewrite  \eqref{ricci100} in the form, 

\begin{eqnarray}
\rr_{\mu\nu}(H)&=&I_{\mu\nu}+II_{\mu\nu}+III_{\mu\nu}+\err\label{ricci102}\\
I_{\mu\nu}&=&\frac{1}{2} h^{\a\b}H_{[\a\b\mu\nu]}\label{I44.1}\\
II_{\mu\nu}&=& [G^{\a\b},P]\,G_{[\a\b\mu\nu]}\\
III_{\mu\nu}&=&- H^{\a\b} H_{\ga\de}
\bigg(\{\Ga^\ga_{\mu\b}\,,\,\Ga^\de_{\a\nu}\} - 
\{\Ga^\ga_{\mu\nu}\,,\,\Ga^\de_{\a\b}\}
\bigg)\\
\err&=&\pi(\pr H,\pr H,\pr G) + \pi(h,\pr G,\pr G)+\pi(h,h;\pr^2 H)
\label{error99}
\end{eqnarray}
\begin{remark}
\label{similarpi} Recalling the definition of $\pi$ and
using the fact that the frequency range  of $h$ is included
in $|\xi|\ge 1$ we have
 \beaa
\int_{u+1}^t\|\pi(h,\pr G,\pr G)\|_{L^2(\dtau)}&\les &
\|h\|_{L_t^2L_x^\infty} \cdot\|\pr G\|_{L_t^2L_x^\infty} 
\cdot \sup_{\tau}\|\pr G\|_{L^2(\dtau)}\\
&\les&\|\pr h\|_{L_t^2L_x^\infty} \cdot\|\pr G\|_{L_t^2L_x^\infty} 
\cdot \sup_{\tau}\|\pr G\|_{L^2(\dtau)}\\
&\les&\|\pr G\|_{L_t^2L_x^\infty} \cdot\|\pr G\|_{L_t^2L_x^\infty} 
\cdot \sup_{\tau}\|\pr G\|_{L^2(\dtau)}
\eeaa
We can thus replace  $\pi(h,\pr G,\pr G)$ by $\pi(\pr G, \pr G; \pr G)$.
By a similar argument, taking into account
the frequency support of $H$,   we can  also replace  $\pi(h,h;\pr^2H)$
by $\pi(\pr G, \pr G; \pr G)$.
 Finally, by a trivial argument, we can also replace  $\pi(\pr H,\pr H,\pr G)$
by $\pi(\pr G, \pr G; \pr G)$. Therefore the error term in \eqref{error99}
can be simplified to
$$\err=\pi(\pr G, \pr G; \pr G).$$
 \end{remark}

\section{The structure of $I_{44}$ and $II_{44}$}
\subsection{Structure of $I_{44}$}
 Contracting the formula \eqref{I44.1} with $L^\mu L^\nu$ 
we obtain, see also  the definition of $[[\,\,]]$ in \eqref{ricci99},
\be{ricci103}
 I_{44}=\frac{1}{2}L^\mu L^\nu \bigg (h^{\a\b}H_{[[\a\b\mu\nu]]} -
 h^{\a\b}H_{\mu\nu\,,\a\b}\bigg)
\end{equation}
Observe  that 
$$L^\mu L^\nu\, H_{[[\a\b\mu\nu]]}= L^\mu L^\nu\,\bigg(H_{\a\nu\,,
\b\mu}+H_{\b\mu\,,
\a\nu}-H_{\a\b\,, \mu\nu}\bigg)\approx L(\pr H)$$

It remains to consider the term
$L^\mu L^\nu \,h^{\a\b}H_{\mu\nu\,,\a\b}$.
 Observe that,
$$
L^\mu L^\nu \pr_\a\pr_\b H_{\mu\nu} = 
\DC \pr H +  [\pr G\pr G]+ [\pr G][\pr G] + [\nab L] \pr G + [\nab L \pr G].
$$
This is obvious if $\a=1,2, 4$ and follows from \eqref{2mod14} of Lemma \ref{Wcord}
 if $\a=3$. Therefore,
$$
L^\mu L^\nu 
 h^{\a\b}H_{\mu\nu\,,\a\b}=
h \DC \pr H  + \pi(h,\pr G;\pr G) + \pi(h,\pr G,\nab L)
$$
Appealing to remark \ref{similarpi} we can 
 summarize our results above in the following 
\begin{proposition} We can write,
$$
I_{44}= h \DC \pr H + \err ,
$$
where 
$$
\err =  \pi(\pr G,\pr G;\pr G) + \pi(\pr G,\pr G;\nab L)
$$
\label{I}
\end{proposition}

\subsection{ Structure of $II_{44}$}\,\,\,
Recall that 
$$
II_{\mu\nu} = \{G^{\a\b}, G_{[\a\b\mu\nu]}\}' = [G^{\a\b}, P] G_{[\a\b\mu\nu]}
$$
For technical reasons we also introduce the following, 
\begin{definition}
\label{fcircw}
Given scalar functions $f,v,w$ we define,
\be{fcw1}
 \{v\,, \,f\circ w\}' = [v, P] \,f\, \bP w +  
\sum_{\la_1>1,\la_2>2, |\ln (\la_1\la_2^{-1})|\le 2} 
[v^{\la_1}, P]\, f\, w^{\la_2} 
\end{equation}
\end{definition}
\begin{lemma}
\label{fcircpr}
\be{ricci5.2}
f \{v\,, \,w\}'= \{v\,, \,f\circ w\}' + \pi(\nab f,\nab v; w)
\end{equation}
\end{lemma}
\begin{proof}
Using representation \eqref{PGPG} for $\{\,,\,\}'$, the commutation lemma 
\ref{Commutation}
and the definition of $\pi$ ( see definition \ref{ric2.989} )
we infer that 
\bea
f \{v\,, \,w\}'& =&  [v, P]\, f\, \bP w +  \sum_{\la_1>1,\la_2>2, |\ln (\la_1\la_2^{-1})|\le 2} 
[v^{\la_1}, P]\, f\, w^{\la_2}\nn \\ &+&  \big[f,[v, P]\big] \bP w + 
 \sum_{\la_1>1,\la_2>2, |\ln (\la_1\la_2^{-1})|\le 2} 
\big[f,\,[v^{\la_1}, P]\big] w^{\la_2}\label{ricci5.1} \\ &=&
\{v\,,\,f\circ w\}' + \pi(\nab f,\nab v;w)\nn
\eea
\end{proof}
We also define,
\beaa
\{v, L\circ w\}'&:=&\{v, L^\mu\circ \pr_\mu w\}'\\
\{v,e_A\circ w\}'&:=&\{v,e_A^\mu\circ \pr_\mu w\}'
\eeaa
\begin{definition} 
\label{fcdc}
We denote by 
$\{v,\DC \circ \,w\}'$ a scalar quantity which can be estimated
as follows
$$|\{v,\DC \circ w\}'|\les 
|\{v, L\circ w\}'|+\big(\sum_{A=1,2}|\{v, e_A\circ w\}'|^2\big)^{\f12}
$$
\end{definition}

We now proceed with the estimate for $II_{44}$.
\beaa
II_{44}&=&L^\mu L^\nu  \{G^{\a\b}, G_{[\a\b\mu\nu]}\}'
\\ &=&L^\mu L^\nu \bigg(\{G^{\a\b}, G_{[[\a\b\mu\nu]]}\}'-
  \{G^{\a\b}, G_{\mu\nu,\,\a\b}\}'\bigg)
\eeaa
We start again with the term containing $ G_{[[\a\b\mu\nu]]}$. According to 
the definition \ref{fcircw}  of $\{v\,,\,f\circ w\}'$ and the relation \eqref{ricci5.2},
 we obtain 
$$
L^\mu L^\nu  \{G^{\a\b}, G_{[[\a\b\mu\nu]]}\}'
=\{G^{\a\b}\,,\, L^\mu L^\nu \circ\, G_{[[\a\b\mu\nu]]}\}' + 
\pi (\nab L, \pr G; \pr^2 G)
$$
Using also  definition \eqref{fcdc},
we infer that 
$$
L^\mu L^\nu \{G^{\a\b}, G_{[[\a\b\mu\nu]]}\}'
=  \{G\,,\, \DC \,\circ\, \pr G\}' + 
\pi(\nab L,\pr G; \pr^2 G)
$$
It remains to consider $L^\mu L^\nu \{G^{\a\b}, G_{\mu\nu,\,\a\b}\}'$. 
Proceeding as above we obtain
$$
L^\mu L^\nu \{G^{\a\b},G_{\mu\nu\,,\a\b}\}'=
\{G^{\a\b}\,,\,  L^\mu L^\nu \circ G_{\mu\nu\,,\a\b}\}'
+ \pi(\nab L, \pr G, \pr^2 G).
$$
According to the wave coordinate condition \eqref{2mod14},
\beaa
L^\mu L^\nu \pr_\a\pr_\b (Q G_{\mu\nu}) &= &
H \cdot \DC \pr (QG)\\
&+&\frac{1}{|Q|}\bigg([\pr G]\c\pr^2 G+[\pr^2 G]\c\pr G\bigg)
+h\c\pr^2 (QG) + \pr G\c\pr (QG)
\eeaa
for any projection 
$Q=I, P, P_{\la_1}$ with $\la_1>1$. Therefore, in view of  definition \eqref{fcw1},
\bea
\{G^{\a\b}\,,\,  L^\mu L^\nu \circ\,\pr_\a\pr_\b G_{\mu\nu}\}'&=&
\{G\,,\,H\c\DC \pr \circ G\}'\label{princ}\\
&+& \{G\,,\, h \circ\pr^2 G\}'   +\{G\,,\, \pr G \circ\pr
G\}'+E \label{ric6.1}
\eea
where the error term $E$ has the form,
\beaa
E&=&    [G, P] \,\, \bP\bigg([\pr G]\pr^2 G  +  [\pr^2 G]\pr G\bigg)\\
&+&\sum_{\la_1>1,\la_2>2, |\ln (\la_1\la_2^{-1})|\le 2} \frac{1}{\la_2}
[ P_{\la_1}G, P]\, \, P_{\la_2}\bigg([\pr G]\pr^2 G  +  [\pr^2 G]\pr G\bigg) 
\eeaa
Observe that the infinite sum above 
is controlled by the presence of the factor
$\la_2^{-1}$ and therefore  $E$ is of the form 
$$E=\pi(\pr G,\pr G\,;\, \pr^2 G).$$
Observe also that the error terms in \eqref{ric6.1} can also be written
 in the form,
\beaa
\{G\,,\, h \circ\pr^2 G\}' &=&\pi(\pr G, h\,;\,  \pr^2 G) \\
 \{G\,,\, \pr G \circ\pr
G\}'&=&\pi(\pr G,\pr G\,;\, \pr G)
\eeaa
Finally, according to lemma \ref{fcircpr} the principal term in \eqref{princ}
$$
\{G\,,\,H\c\DC \pr \circ G\}' = H\c \{G\,,\,\DC \pr \circ G\}' +
\pi(\nab H, \pr G, pr^2 G)
$$
We summarize these calculations in the following.
\begin{proposition}
We can write
\be{eq:II44}
II_{44} = \{G\,,\, \DC \pr \circ\, G\}' + H\c \{G\,,\,\DC \pr \circ G\}'+ \err,
\end{equation}
where the error term
$$
\err  = \pi(\nab L,\pr G; \pr^2 G) + 
\pi(\pr G,\pr G; \pr G) +
\pi(\pr G, \pr G; \pr^2 G)
$$
\label{II}
\end{proposition}

\section{The structure of $III_{44}$}
 Recall  \eqref{2mod1},
\bea
 III_{44}&=&-  L^\mu L^\nu 
\,H^{\a\b} H_{\ga\de}
\bigg(\big\{\Ga^\ga_{\mu\b}(G), \Ga^\de_{\a\nu}(G)\big\} - 
\big\{\Ga^\ga_{\mu\nu}(G),\Ga^\de_{\a\b}(G)\big\}\bigg)\nn\\
&=&\EE_1-\EE_2\label{2mod15}
\eea
with $\{\,\,,\,\}$ denoting the modified paradifferential product
introduced in definition \ref{paradifferential}. 
\begin{remark} We note here the following
simple property of $\{\,\,,\,\}$:
$$
f\{v,w\}= \{fv,w\} + \pi(v, w; \nab f)=  \{v,f w\} + \pi(v,w;\nab f).
$$
\label{remark7.1}
\end{remark}
Recalling remark \ref{remark5.5} we shall now introduce
the following expression closely related to $fg\{v,w\}$.
\begin{definition}  Given scalars $v,w,f,g$ we introduce
\begin{align}
\{f\circ v,g\circ w\}&= P\,\sum_{\half<\la_1\le 4 }
fv^{\la_1}\c gP_{\le \f12}\, w
+ P\,\sum_{\la_1>\half, |\ln (\la_1\la_2^{-1})|\le 2}
fv^{\la_1}\,\c  gw^{\la_2}\nn  
 \\ &+ P\,\sum_{\half < \la_2\le 4 }
 fP_{\le \f12} v\,\c gw^{\la_2} 
+ P\,\sum_{\la_2 >\half, |\ln (\la_1\la_2^{-1})|\le 2}
fv^{\la_1}\,\c gw^{\la_2}\label{figdec'}\\
& - \sum_{\half<\la_1\le 1 }
fv^{\la_1}\,\c gP_{\le \f12}w
- \sum_{\half <\la_2\le 1}
fP_{\le \f12}v\,\c gw^{\la_2}\nn  
 \\ &- \sum_{\half < \la_1, \la_2\le 1 }
fv^{\la_1} \,\c gw^{\la_2}\nn
\end{align}
\label{figure}
\end{definition}
\begin{lemma} We have,
\be{ric7.99}
f\{v\,,\, w\}=\{f\circ v\,,\,w\}+\pi( v, w;\nab f)=
\{ v\,,\,f\circ w\}+\pi( v, w;\nab f)
\end{equation}
\label{fwbrac}
\end{lemma}
\begin{proof}
\end{proof}
We also define,
\beaa
\{L\circ v,w\}&:=&\{L^\mu\circ \pr_\mu v, w\}\\
\{e_A\circ v,w\}&:=&\{e_A^\mu\circ \pr_\mu v,w\}
\eeaa
In view of the lemma we have 
\be{ric100}
L^\mu\{\pr_\mu v,w\}=\{L\circ v,w\}+\pi(\pr v,w; \nab L)
\end{equation}
\begin{definition} We denote by 
$\{\DC \circ v,w\}$ a scalar quantity which can be estimated
as follows
$$|\{\DC \circ v,w\}|\les |\{L\circ v,w\}|+\big(\sum_{A=1,2}|\{e_A\circ
v,w\}|^2\big)^{\f12}
$$
\end{definition}

In the calculation below we shall
use the notation  $\Ga^{\ga}_{\a\b}= G^{\ga\si}\Ga_{\si|\a\b}$
where,
$$\Ga_{\si|\a\b}=\f12\big(G_{\a\si,\,\b}+G_{\b\si,\,\a}-G_{\a\b,\,\si}\big).$$

\medn
{\em The  term $\EE_1=-H^{\a\b} H_{\ga\de}L^\mu L^\nu 
\bigg\{\Ga^\ga_{\mu\b}(G), 
\Ga^\de_{\a\nu}(G)\bigg\}$:}\,\,\,\,

Using remark \ref{remark7.1}  then  expressing $G^{\a\b}=H^{\a\b}-h^{\a\b}+O(h^2)$
 and applying the definition of $\pi$ we derive,
\bea
-\EE_1&:=&H^{\a\b} H_{\ga\de}L^\mu L^\nu  \bigg\{\Ga^\ga_{\mu\b}(G), 
\Ga^\de_{\a\nu}(G)\bigg\}\label{ricci7.51}\\ &=& H^{\a\b} H_{\ga\de}L^\mu L^\nu  
\bigg\{ G^{\ga\rho}\Ga_{\rho|\mu\b}\,, \,G^{\de\si}\Ga_{\si|\a\nu}   \bigg\}    \nn \\ 
&=& H^{\a\b} H_{\ga\de} G^{\ga\rho}G^{\de\si}L^\mu L^\nu  
\bigg\{\Ga_{\rho|\mu\b}\,, \,\Ga_{\si|\a\nu}   \bigg\}+\err\nn\\
&=&H^{\a\b}H^{\de\si}L^\mu L^\nu\bigg\{\Ga_{\de |\mu\b}\,, \,\Ga_{\si|\a\nu}  
\bigg\}+\err
\nn\\ &=&H^{\a\b}H^{\de\si}\bigg \{L^\mu\circ\Ga_{\de |\mu\b}\,,
\,L^\nu\circ \Ga_{\si|\a\nu}  
\bigg\}+\err\quad \mbox{(using \eqref{ric7.99})}\nn
\eea
with the final expression\footnote{Use also remark \ref{similarpi}, lemma 
\ref{fwbrac},
and the boundedness of $\|G, H, h\|_{L^\infty}$.} for  the error term
$$\err=\pi(\pr G,\pr G\,;\,\pr G) + \pi(\pr G, \pr G\,;\, \nab L) $$

Consider now the bilinear  term $\{L^\mu\circ\Ga_{\de |\mu\b}\,,
\,L^\nu\circ \Ga_{\si|\a\nu}\} $. As we start manipulating 
the left hand side we consider $\{L^\mu\circ\Ga_{\de |\mu\b}\,,w\}$ for a fixed
$w$. As $w$ remains unchanged in the calculations below we shall drop the 
bracket and  simply write
$\{L^\mu\circ\Ga_{\de |\mu\b}\,,w\}=L^\mu\circ\Ga_{\de |\mu\b}$. 
Thus instead of,
$$\bigg\{L^\mu\circ \Ga_{\de |\mu\b}\,,\,w\bigg\}=\f12\bigg\{ L^\mu \circ
\bigg(G_{\mu\de,\,\b} + G_{\b\de,\,\mu}- G_{\b\mu,\,\de}\bigg)\,,\, w\bigg\}=\cdots$$
we write,
\beaa
L^\mu\circ \Ga_{\de |\mu\b}&=&\f12 L^\mu \circ \bigg(G_{\mu\de,\,\b} + G_{\b\de,\,\mu}-
G_{\b\mu,\,\de}\bigg)=\f12 L^\mu \circ \bigg(G_{\mu\de,\,\b} -
G_{\b\mu,\,\de}\bigg)+\DC\circ G\\ &=&
\f12 \bigg( L^\mu\pr_\b\circ  G_{\mu\de} -
L^\mu \pr_\de\circ G_{\b\mu}\bigg)+\DC\circ G,
\eeaa
where we have used that $\pr$ commutes with $\circ$, i.e. 
$\{f\circ \pr v\,,\, w\}=\{f\pr \circ v\,,\, w\}$  .
Recall that, see \eqref{decompartial}, $\pr_\a=-\f12 L_\a \, \Lb +\DC.$
Therefore,
\be{ricci7.1}
L^\mu\circ \Ga_{\de |\mu\b}=-\frac 14 \big ({L}_{\b} L^\mu\pr_3 \circ G_{\mu\de} 
- {L}_\de L^\mu\pr_3\circ G_{\b\mu}\bigg) + \DC\circ G 
\end{equation}

According to \eqref{2mod11} of Lemma \ref{Wcord} and
the formula $H^{3\a}=-\f12 L^\a, $ we have
$$
- L^\a \pr_3 (QG)_{\a\si}= H^{\a\b} \pr_\si (QG)_{\a\b} +G\DC(Q G) +
h \pr(Q G)  + \frac{1}{|Q|}[\pr G]\pr G
$$
 with $Q$  any  of the projections $Q=I, P, P_{\la_1}$, with $\la_1>1$,
 appearing in the definition of $\,\{\,\,,\,\,\}$ and $\circ$. Therefore,
$$
- L^\a \pr_3 \circ G_{\a\si}= H^{\a\b} \pr_\si \circ G_{\a\b} +G\DC\circ G +
h\circ \pr G  + [\pr G]\pr G
$$
Therefore, from \eqref{ricci7.1},  and expanding $\pr_\de, \pr_\b$
relative to the null frame,
$$
\aligned
 L^\mu\circ \Ga_{\de |\mu\b}& = \frac 14 \big (L_{\b}H^{\mu\si}\pr_\de
\circ G_{\mu\si} - L_\de H^{\mu\si}\pr_\b \circ G_{\mu\si}\big) + 
  \DC \circ G \\ &+  G\cdot \DC \circ G +  h\circ \pr G + [\pr G]\pr G\\ &=
-\frac 18 \big (L_\b L_\de H^{\mu\si}\pr_3 \circ G_{\mu\si} -
L_\de L_\b H^{\mu\si}\pr_3 \circ G_{\mu\si}\big)\\ & + 
 \DC \circ G +  G\cdot \DC\circ  G + h\circ \pr G +[\pr G]\pr G\\ & = 
\DC \circ G +  G\cdot \DC\circ  G + h\circ \pr G +[\pr G]\pr G 
\endaligned
$$
Similarly we  have
\beaa
L^\nu \circ \Ga_{\a|\nu\ga}&=&\f12
L^\nu\circ \bigg(G_{\nu\ga,\,\a} + G_{\a\ga,\,\nu}-
G_{\a\nu,\,\ga}\bigg)\\
 &=& -\frac 14\bigg( L^\nu {L}_{\a}\pr_3 \circ  G_{\nu\ga} -
L^\nu {L}_{\ga}\pr_3 \circ  G_{\a\nu}\bigg) + \DC\circ G \\ &=& 
\frac 14 \bigg( {L}_{\a} H^{\nu\si} \pr_\ga \circ G_{\nu\si} -
{L}_{\ga} H^{\nu\si}\pr_\a \circ G_{\nu\si}\bigg)\\ & +& \DC\circ G + G\cdot \DC\circ G +
 h\circ \pr G + [\pr G]\pr G  
\\ &=&  \DC\circ G + G\cdot \DC\circ G +
 h\circ \pr G + [\pr G]\pr G
\eeaa
Thus, going back to \eqref{ricci7.51},
\begin{align}
\EE_1
&=  H\cdot H \bigg\{(\DC\circ G + G\cdot \DC\circ G),(\DC\circ G + 
G\cdot \DC\circ G)\bigg \}\nn\\
&+  \pi(\pr G,\pr G\,;\,\pr G) + \pi(\pr G,\pr G\,;\,\nab L)\label{mod161} 
\end{align}

\medn
{\em The term $\EE_2=H^{\a\b} H_{\ga\de}L^\mu L^\nu \bigg\{\Ga^\ga_{\mu\nu}(G),
\Ga^\de_{\a\b}(G)\bigg\}$:}\,\,\,\,

 Using remarks \ref{remark7.1} and \ref{fwbrac},
\beaa
\EE_2
&=&H^{\a\b} H_{\ga\de} L^\mu L^\nu 
\bigg\{H^{\ga \eps}\Ga_{\ep|\mu\nu}\,,\,H^{\de\si}\Ga_{\si|\a\b}\bigg\}+
\pi(h, \pr G\,;\pr G)\,\\ &=& H^{\de\si} \bigg\{ L^\mu L^\nu\circ \Ga_{\de|\mu\nu} \,   ,  \,
H^{\a\b}\circ \Ga_{\si|\a\b}  \bigg\}+\pi(\pr G,\pr G; \pr G) + \pi(\pr G, \pr G; \nab L)
\eeaa

Observe that according to \eqref{2mod1000}, Lemma \ref{Wcord},
$$
 L^\mu L^\nu \pr_\de(QG)_{\mu\nu} = G\cdot \DC (QG) + h\pr (QG) +
\frac{1}{|Q|}[\pr G]\pr G
$$
for any $Q=I, P, P_{\la_1}$ with $\la_1>1$.
Therefore,
$$
 L^\mu L^\nu \pr_\de\circ G_{\mu\nu} = G\cdot \DC \circ G + h\circ \pr G +[\pr G]\pr G
$$
Thus 
\bea
 L^\mu L^\nu\circ \Ga_{\de|\mu\nu} &=&
 L^\mu L^\nu \circ \bigg(G_{\mu\de,\,\nu} + 
G_{\nu\de,\,\mu} - G_{\mu\nu,\,\de}\bigg)\nn\\
&= &\DC \circ G -
 L^\mu L^\nu \pr_\de \circ G_{\mu\nu} = \DC\circ  G +  G\cdot \DC\circ 
G\label{ricci7.7}\\ & +& 
 h\circ \pr G +[\pr G]\pr G \nn
\eea
Using \eqref{2mod1} we  have,
$$H^{\a\b}\pr_\si \circ G_{\a\b}=2H^{3\a}\pr_3 \circ G_{\a\si}+G\DC \circ G+
h\circ\pr G +[\pr G]\pr G.$$
Therefore,
\beaa
H^{\a\b}\circ \Ga_{\si|\a\b}&=&H^{\a\b}\circ \bigg( G_{\a\si,\,\b} + G_{\b\si,\,\a} -
G_{\a\b,\,\si}\bigg)\\
& = & H^{\a\b} \circ \bigg( 2 G_{\a\si,\,\b} - G_{\a\b,\,\si}\bigg)\\
&=&2 H^{\a\b}\pr_\b\circ  G_{\a\si}- 2H^{3\a}\pr_3 \circ G_{\a\si}+G\DC \circ G +
h\circ \pr G +[\pr G]\pr G\\
 &=& G\DC\circ  G +h\circ \pr G +[\pr G]\pr G
\eeaa
Therefore, similar to \eqref{mod161}, we derive
\bea
\EE_2&=& H \bigg\{(\DC\circ G + G\cdot \DC \circ G)\,,\,(\DC\circ G + 
G\cdot \DC\circ G)\bigg \}\label{mod162}\\
&+&  \pi(\pr G,\pr G\,;\,\pr G) + \pi(\pr G,\pr G; \nab L) \nn
\eea
We now observe that according to the remark \ref{remark7.1} 
$$
\{G\cdot \DC\circ G, f\}= G \{\DC\circ G, f\} + \pi(\pr G, f; \nab G)
$$ 
Therefore returning to \eqref{2mod15}, using \eqref{mod161}, 
\eqref{mod162}, and the boundedness of $H$ and $G$  we infer the following
\begin{proposition}
We can write
$$
 III_{44} =  \,\{\DC\circ G\,,\,\DC\circ G\} +\err,
$$
where
$$
\err =  \pi(\pr G,\pr G\,;\, \pr G) + \pi(\pr G,\pr G\,;\,\nab L) 
$$
\label{III}
\end{proposition}
\section{Estimates for $I_{44}$, $II_{44}$, and $III_{44}$}
According to the reduction \eqref{R0} and the representation
\eqref{ricci102},
$$
\rr_{44} = I_{44} +II_{44} +III_{44} + \err
$$
with 
$$
\err =\pi(\pr H,\pr H;\pr G) + \pi(h,\pr G;\pr G)+\pi(h,h;\pr^2 H).
$$
Therefore we need to show that 
\beaa
\int_{u+1}^t\|I_{44}\|_{L^2(\dtau)}\,d\tau &+& 
 \int_{u+1}^t\|II_{44}\|_{L^2(\dtau)}\,d\tau\\ &+&
\int_{u+1}^t\|III_{44}\|_{L^2(\dtau)}\,d\tau +
\int_{u+1}^t\|\err\|_{L^2(\dtau)}\,d\tau \les \la^{-1}
\eeaa
We start with error terms accumulated above and in the lemmas \ref{I}, \ref{II},
\ref{III}.
\subsection{Estimates for the error terms}\,\,\,
According to the property \eqref{tric3.4} of $\pi$, 
$$
\pi(\pr H,\pr H;\pr G)\le \pi(\pr G,\pr G;\pr G).
$$ 
We then estimate, with the help of the  estimates  
\eqref{asG1}--\eqref{asG3} for $G$,
\beaa
\int_{u+1}^t\|\pi(\pr G,\pr G;\pr G)\|_{L^2(\dtau)}\,d\tau &\les&
\|\pr G\|^2_{L^2_t L^\infty_x} \sup_{\tau} \|\pr G\|_{L^2(\dtau)}\\
&\les& \la^{-1-8\eps_0} \sup_{\tau, u} \|\pr G\|_{L^2(\stau)} \les\la^{-1-10\eps_0}
\eeaa
Since the frequencies of $h$ are restricted to the region 
$|\xi|\ge 1$, $h=(I-P) G$, we also have 
\beaa
\int_{u+1}^t\|\pi(h,\pr G;\pr G)\|_{L^2(\dtau)}\,d\tau &\les& 
\|h\|_{L^2_t L^\infty_x} \|\pr G\|_{L^2_t L^\infty_x} \sup_{\tau} \|\pr G\|_{L^2(\dtau)}
\\ &\les& \la^{-\frac 12-6\eps_0} \|\pr h\|^2_{L^2_t L^\infty_x}\les 
\la^{-1-10\eps_0}
\eeaa
In addition, using the background estimates \eqref{aso1}--\eqref{aso3},
\beaa
\int_{u+1}^t\|\pi(h,h;\pr^2 H)\|_{L^2(\dtau)}\,d\tau &\les& 
\|h\|^2_{L^2_t L^\infty_x}\sup_{\tau} \|\pr^2 H\|_{L^2(\dtau)}\\
&\les & \|\pr h\|^2_{L^2_t L^\infty_x}\sup_{\tau} \|\pr^2 H\|_{L^2(\Si_{\tau})}\les
\la^{-\frac 32-8\eps_0}
\eeaa
Estimating the error terms generated in proposition \ref{I}, 
and using the  estimate \eqref{DL} for   $\nab L$
\beaa
\int_{u+1}^t\|\pi(\pr G,\pr G;\nab L)\|_{L^2(\dtau)}\,d\tau &\les& 
\|\pr G\|^2_{L^2_t L^\infty_x} \sup_{\tau} \|\nab L\|_{L^2(\dtau)}\les
\la^{-1-8\eps_0}
\eeaa
To bound the error term $\pi(\nab L, \pr G; \pr^2 G)$ in  proposition \ref{II}
we use the inequality  \eqref{nabL}, $|\nab L|\les (\Theta + r^{-1})$ and
\beaa
\bigg(\int_{u+1}^t \|(\Theta + r^{-1})\|^2_{L^\infty(\dtau)}\,d\tau\bigg)^{\frac 12}
&\les &\|\Theta\|_{L^2_t L^\infty_x} + 
\bigg(\int_{u+1}^t \frac {d\tau}{(\tau -u)^2}\bigg)^{\frac 12}
\\ &\les& \la^{-\frac 12 -2\eps_0} +1, 
\eeaa
which follows from the comparison $r\approx \tau -u$, see 
\eqref{ric2.990}.
Thus,
\beaa
\int_{u+1}^t\|\pi(\nab L,\pr G;\pr G)\|_{L^2(\dtau)}\,d\tau &\les& 
\bigg(\int_{u+1}^t \|\nab L\|^2_{L^\infty(\dtau)}\,d\tau\bigg)^{\frac 12}
\|\pr G\|_{L^2_t L^\infty_x} \sup_{\tau} \|\pr^2 G\|_{L^2(\dtau)}\\&\les&
\la^{-\frac 12 -4\eps_0} \sup_{\tau} \|\pr^2 G\|_{L^2(\Si_\tau)}
\les  \la^{-1-4\eps_0}
\eeaa
Finally,
\beaa
\int_{u+1}^t\|\pi(\pr G,\pr G;\pr^2 G)\|_{L^2(\dtau)}\,d\tau &\les& 
\|\pr G\|^2_{L^2_t L^\infty_x}\sup_{\tau} \|\pr^2 G\|_{L^2(\dtau)}\\
&\les & \la^{-1-8\eps_0} \sup_{\tau} \|\pr^2 G\|_{L^2(\Si_\tau)}
\les \la^{-\frac 32 -8\eps_0}
\eeaa
The error terms in  proposition \ref{III} are the same as considered above.
\subsection{Estimates for the principal terms}
These estimates depend decisively on the $L^2(C_u)$ estimates
for the  tangential  derivatives of  $G$ and $H$ derived in 
proposition 7.7 of \cite{Einst2}, see also proposition   \ref{CuH}. 
For convenience we
recall the result here.
\be{cuH}
\|D_*\pr H\|_{L^2(C_u)} \les \la^{-\f12},\qquad
\|D_* H\|_{L^2(C_u)} \les \la^{\f12}
\end{equation}
Also,
\bea
\|D_*\pr (P_\mu G)\|_{L^2(C_u)} &\les& \mu^{\frac 12 - 4\eps_0} \la^{-\frac 12 -
4\eps_0},\nn\\ \|D_* (P_{\mu}G)\|_{L^2(C_u)}&\les& \la^{-\frac 12-4\eps_0}
\mu^{-\frac 12-4\eps_0}\label{10cuH}
\eea

We start with the principal term $h\DC \pr H$ 
appearing in proposition \ref{I}.
\beaa
\int_{u+1}^t\|h\DC \pr H\|_{L^2(\dtau)}\,d\tau &\les&
\int_{u+1}^t\|h(\tau)\|_{L^\infty_x}\|\DC \pr
H\|_{L^2(\dtau)} d \tau\\
&\les &\|h\|_{L_t^2L_x^\infty}\bigg(\int_{u+1}^t
\|\DC \pr H\|_{L^2(\dtau)}^2\,d\tau\bigg)^{\f12}\\
&\les &\|\pr G\|_{L_t^2L_x^\infty}\sup_{u\le u'\le u+1}
\|\DC \pr H\|_{L^2(C_{u'})}\\
&\les&  \la^{-1-4\eps_0}\qquad \qquad \mbox{(using \eqref{cuH})}
\eeaa
as desired. 

We now estimate the principal terms 
$\{G\,,\, \DC \pr \circ\, G\}'$  and $H\c \{G\,,\,\DC \pr \circ G\}'$
appearing in proposition \ref{II}. Since $H$ is bounded 
it clearly suffices to treat the first term.
Recall that 
$$
\{G\,, \, \DC \pr\circ  G\}' = [G, P] \,\DC\, \pr(\bP G) +  
\sum_{\nu>1,\mu>2, |\ln (\nu\mu^{-1})|\le 2} 
[P_{\nu}G, P]\, \DC\,\pr (P_{\mu} G) 
$$
We estimate the first term as follows: 
\bea
\int_{u+1}^t\|[G, P] \,\DC\, \pr(\bP G)\|_{L^2(\dtau)}\,d\tau &\les& 
\int_{u+1}^t\|\pr G\|_{L^\infty_x} \|\DC\, \pr(\bP G)\|_{L^2(\dtau)}\,d\tau\nn \\
&\les & \|\pr G\|_{L^2_t L^\infty_x} 
\bigg(\int_{u+1}^t\|\DC\, \pr(\bP G)\|^2_{L^2(\dtau)}\bigg)^{\frac 12}\nn\\
&\les &\|\pr G\|_{L_t^2L_x^\infty}\sup_{u\le u'\le u+1}
\|\DC \pr (\bP G) \|_{L^2(C_{u'})}\nn \\ &\les&  \la^{-1-4\eps_0} \qquad \mbox{(using
\eqref{cuH})}
\label{ric9.1} 
\eea
 We estimate the high-high interaction as follows
\beaa
\int_{u+1}^t\|[P_{\nu}G, P] \,\DC\, \pr(P_{\mu} G)\|_{L^2(\dtau)}\,d\tau &\les& 
\int_{u+1}^t \|P_\nu G\|_{L_x^\infty}\|\DC\, \pr(P_\mu G)\|_{L^2(\dtau)}\\
&\les & \frac 1\nu 
\int_{u+1}^t \|P_\nu \nab G\|_{L_x^\infty}\|\DC\, \pr(P_\mu G)\|_{L^2(\dtau)}\\
&\les &  \frac 1\nu \|\pr G\|_{L_t^2L_x^\infty}\sup_{u\le u'\le u+1}
\|\DC \pr (P_\mu G) \|_{L^2(C_{u'})}\\
 &\les& 
\nu^{-1} \mu^{\frac 12 - 4\eps_0} \la^{-1 - 8\eps_0}\qquad \mbox{(using
\eqref{10cuH})}
\eeaa

Therefore,
\be{ric9.2}
\int_{u+1}^t\|[P_{\nu}G, P] \,\DC\, \pr(P_{\mu} G)\|_{L^2(\dtau)}\,d\tau \les 
\nu^{-1} \mu^{\frac 12 - 4\eps_0} \la^{-1 - 8\eps_0}
\end{equation}
Combining \eqref{ric9.1} and \eqref{ric9.2} we conclude that 
\beaa
\int_{u+1}^t \|\{G\,, \, \DC \pr\,\circ\,  G\}' \|_{L^2(\dtau)}\,d\tau&\les& 
 \la^{-1-4\eps_0} + \,\,\, \la^{-1-8\eps_0} 
\sum_{\nu>1,\mu>2, |\ln (\nu\mu^{-1})|\le 2} \nu^{-1} \mu^{\frac 12 - 4\eps_0}\\
&\les &  \la^{-1-4\eps_0}
\eeaa

It remains to estimate the principal term  $\{\DC\circ G\,,\, \DC\circ G\}$ in 
Proposition \ref{III}. We recall from definition \ref{figure} that 
\beaa
\{\DC\circ G\,,\, \DC\circ G\}&=& P\,\sum_{\half<\nu\le 4 }
\DC (P_\nu G) \c \DC (P_{\le \f12}\, G)\\
&+& P\,\sum_{\nu>\half, |\ln (\nu\mu^{-1})|\le 2}
\DC (P_\nu G) \,\c  \DC (P_\mu G)\nn  
 \\ &+& P\,\sum_{\half < \mu\le 4 }
 \DC (P_{\le \f12} G)\,\c \DC (P_\mu G)\\
&+& P\,\sum_{\mu_2 >\half, |\ln (\nu\mu^{-1})|\le 2}
\DC (P_\nu G)\,\c \DC (P_\mu G)\\
& -& \sum_{\half<\nu \le 1 }
\DC (P_\nu G)\,\c \DC( P_{\le \f12}G)
\\ &-& \sum_{\half <\mu\le 1}
\DC (P_{\le \f12}G)\,\c \DC (P_\mu G)  
 \\ &-&\sum_{\half < \nu, \mu \le 1 }
\DC (P_\nu G) \,\c \DC (P_\mu G)
\eeaa
By symmetry and similarity it suffices to estimate the first 2 terms in the 
expression above. 
We have
\beaa
\int_{u+1}^t\|\DC (P_\nu G) \c \DC (P_{\le \f12}\, G)\|_{L^2(\dtau)}\,d\tau &\les&
\int_{u+1}^t \|\DC (P_\nu G)\|_{L^2(\dtau)}\| \DC (P_{\le \f12}\, G)\|_{L^\infty_x}
d\tau\\  &\les& \|\pr G\|_{L_t^2L_x^\infty} \sup_{u\le u'\le u+1}
\|\DC (P_\nu G) \|_{L^2(C_{u'})}\\
&\les&\nu^{-\frac 12-4\eps_0} \la^{-1-8\eps_0} \qquad\qquad \mbox{(by \eqref{10cuH}.)}
\eeaa
Thus
\beaa
\int_{u+1}^t\|P\,\sum_{\half<\nu\le 4 } \DC (P_\nu G) \c \DC (P_{\le \f12}\, G)\|_{L^2(\dtau)}\,
d\tau &\les& \sum_{\half<\nu\le 4 }\int_{u+1}^t
\|\DC (P_\nu G) \c \DC (P_{\le \f12}\, G)\|_{L^2(\dtau)}\\ &\les& \la^{-1-8\eps_0}
\eeaa
Consider now the
 high-high interaction term
$$J=P\,\sum_{\nu>\half, |\ln (\nu\mu^{-1})|\le 2}
\DC (P_\nu G) \c \DC (P_\mu\, G)$$
Clearly,
\beaa
\int_{u+1}^t\|\DC (P_\nu G) \c \DC (P_\mu\, G)\|_{L^2(\dtau)}\,d\tau &\les&
\int_{u+1}^t\|\DC (P_\nu G)\|_{L^2(\dtau)} \| \DC (P_\mu\, G)\|_{L^\infty_x}\,d\tau\\
 &\les& \|\pr G\|_{L_t^2L_x^\infty} \sup_{u\le u'\le u+1}
\|\DC (P_\nu G) \|_{L^2(C_{u'})}\\ &\les& \la^{-1-8\eps_0} \nu^{-\frac 12-4\eps_0}
\eeaa
Thus,
\beaa
\int_{u+1}^t\|\,\,\,J\|_{L^2(\dtau)}\,d\tau &\les&
\sum_{\nu>\half, |\ln (\nu\mu^{-1})|\le 2} 
\int_{u+1}^t\|\DC (P_\nu G) \c \DC (P_\mu\, G)\|_{L^2(\dtau)}\,d\tau\\
 &\les& \la^{-1-8\eps_0} \sum_{\nu>\half, |\ln (\nu\mu^{-1})|\le 2} \nu^{-\frac 12-4\eps_0}
\les \la^{-1-8\eps_0}
\eeaa

\end{document}